\documentclass{article}
\usepackage{graphicx,amsmath,amssymb,amsthm} % Required for inserting images
\usepackage{xcolor,float}

\title{Orbits by the up-down action of braid diagrams}
\author{ Komal Negi, Ayaka Shimizu, Yoshiro Yaguchi, Madeti Prabhakar }
\date{}
\newtheorem{theorem}{Theorem}[section]
\newtheorem{lemma}[theorem]{Lemma}
\newtheorem{prop}[theorem]{Proposition}
\newtheorem{corollary}[theorem]{Corollary}

\theoremstyle{definition}
\newtheorem{definition}[theorem]{Definition}
\newtheorem{remark}[theorem]{Remark}
\newtheorem{Example}[theorem]{Example}
\begin{document}

\maketitle
\begin{abstract}
The set of all virtual or classical braid diagrams forms a monoid and gives a natural monoid action on a direct product of ${\mathbb Z}$ called the up-down action. In this paper, we determine the orbit of every tuple of ${\mathbb Z}$ under the up-down action of virtual or classical braid diagrams. Moreover, we determine the orbit for irreducible braid diagrams. We also consider the isotropy submonoid and give a condition for a braid diagram to admit an up-down coloring to its closure.\\


\end{abstract}

\section{Introduction}

To study knots and links, it is crucial to study colorings for a link diagram which give the arcs elements of a set $S$ according to a certain rule. 
Colorings for braid diagrams are helpful to understand the colorings for link diagrams because each link is a closure of a braid and the set of braids in $\mathbb{R}^3$ or the set of braid diagrams in $\mathbb{R}^2$ has an algebraic structure up to isotopy, a group or a monoid, respectively, and it is suitable to deal with. 

For example, when $S$ is a quandle, a set with a binary operation which axiomatizes the notion of conjugation in a group, the quandle colorings for virtual link diagrams are studied to find link invariants. Quandle colorings for classical links can be studied considering a natural action by a braid group $B_n$ of degree $n$ on the $n$-direct product $S^n$ of $S$, which is called Hurwitz action (see, for example, \cite{D}). 
Although it has been proved that there does not exist any common algorithm for determining all elements in the Hurwitz orbit for general $S$ (\cite{LM}), there are many studies for determination of all the elements in certain orbits by Hurwitz action for specific quandles $S$ (see, for example, \cite{BGRW, Y}). 

On the other hand, the {\it up-down coloring} is a coloring to arcs of a virtual link diagram by elements of $S= \mathbb{Z}$, which is not a quandle coloring, and it was used to determine the necessity of a local move on diagrams (\cite{KAY}). 
For virtual braid diagrams and the twisted virtual braid diagrams, which is a generalization of virtual braids defined in \cite{KPS}, recently the up-down coloring, or up-down labeling, was introduced in \cite{KAM}.  
The up-down coloring for a classical braid diagram (resp. a virtual braid diagram) gives a monoid action on $\mathbb{Z}^n$ by the monoid $BD_n$ (resp. $VBD_n$) formed by classical braid diagrams (resp. virtual braid diagrams) of degree $n$ up to an isotopy of $\mathbb{R}^2$, which is called the {\it up-down action}. We remark here that the up-down action is not a group action. 
In this paper, we determine the orbits $\vec{x}\cdot BD_n$ and $\vec{x}\cdot VBD_n$ of the up-down action by $BD_n$ and $VBD_n$ for any $\vec{x} \in \mathbb{Z}^n$ in the following theorems. 
Let $\mathrm{tr}(\vec{x})=\sum_{i=1}^n x_i$ for $\vec{x}=(x_1, x_2, \dots , x_n) \in \mathbb{Z}^n$, and let $IVBD_n$ (resp. $IBD_n$) be the set of all virtual (resp. classical) braid diagrams of degree $n$ which have no non-alternating bigons. 

\begin{theorem}
For $\vec{x}=(x_1, x_2, \dots , x_n) \in \mathbb{Z}^n$ ($n \geq 3$), 
$$\vec{x}\cdot BD_n = \left\{ \vec{y} \in \mathbb{Z}^n \ | \ \mathrm{tr}(\vec{x})=\mathrm{tr}(\vec{y}), \ I(\vec{x})=I(\vec{y}) \right\} ,$$
where $I(\vec{x})$ is the number of components $x_i$ of $\vec{x}$ such that ``$i$ is odd and $x_i$ is even'' or ``$i$ is even and $x_i$ is odd''. Moreover, $\vec{x}\cdot IBD_n=\vec{x}\cdot BD_n$.
\label{thm-main1}
\end{theorem}

\begin{theorem}
For $\vec{x}=(x_1, x_2, \dots , x_n) \in \mathbb{Z}^n$, 
$$\vec{x}\cdot VBD_n = \left\{ \vec{y} \in \mathbb{Z}^n \ | \ \mathrm{tr}(\vec{x})=\mathrm{tr}(\vec{y}) \right\}.$$ Moreover, $\vec{x}\cdot IVBD_n=\vec{x}\cdot VBD_n$.
\label{thm-main2}
\end{theorem}

The structure of the article is as follows.
In Section 2, we provide essential background information on braid diagrams, the up-down coloring, and the monoid structure of virtual and classical braid diagrams, setting the foundation for our research.
In Section 3, we find the orbit by the up-down action to classical braid diagrams. 
The section is divided into three parts. 
We first try to solve the orbit problem for a submonoid, the set of classical pure braid diagrams $PBD_n$, which helps us in solving the main problem of this section, finding the orbits $\vec{0}\cdot PBD_n$ and $\vec{0}\cdot IPBD_n$ of the zero element $\vec{0}\in {\mathbb Z}^{n}$ under the action of $PBD_n$ and $IPBD_n =IBD_n\cap PBD_n$, respectively. Second, we provide a necessary and sufficient condition for $\vec{y}\in \vec{x}\cdot BD_n$. Third,
we find the orbits $\vec{x}\cdot BD_n$ and $\vec{x}\cdot IBD_n$ of general element $\vec{x}\in {\mathbb Z}^n$ under the action of $BD_n$ and $IBD_n$, respectively. In Section 4, we find the orbits
$\vec{x}\cdot VBD_n$ and $\vec{x}\cdot IVBD_n$ of $\vec{x}\in {\mathbb Z}^n$ 
under the up-down
action of $VBD_n$ and $IVBD_n$, respectively. In Section 5, we study the relation between the OU matrix and up-down action, and consider the isotropy submonoid of $\vec{0}\in {\mathbb Z}^n$ under the up-down action of $VBD_n$. As an application, we give a sufficient condition for an element of $VBD_n$ to admit the up-down coloring to its closure. 
In Section 6, we see the up-down action for some braid classes known as torus and weaving braids.

\section{Preliminaries}

\subsection{Virtual braid diagrams and Classical braid diagrams}
\begin{figure}[ht]
	\centering
\includegraphics[width=9cm]{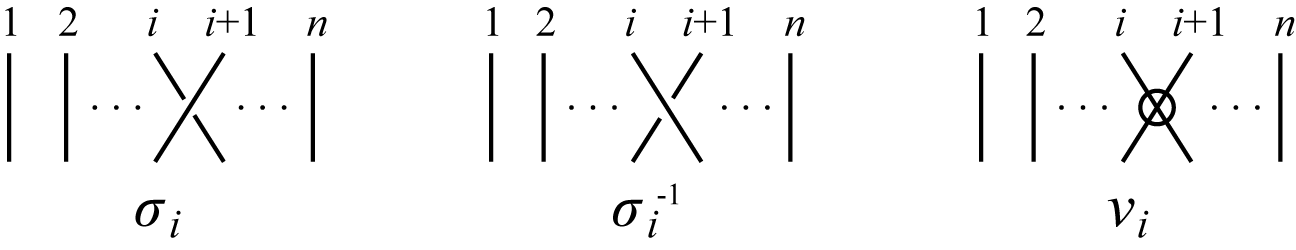}	
  \caption{$\sigma_i$, $\sigma_i^{-1}$ and $v_i$.}
 \label{fig-sigma}
\end{figure}

\begin{figure}[ht]
\centering
\includegraphics[width=8cm]{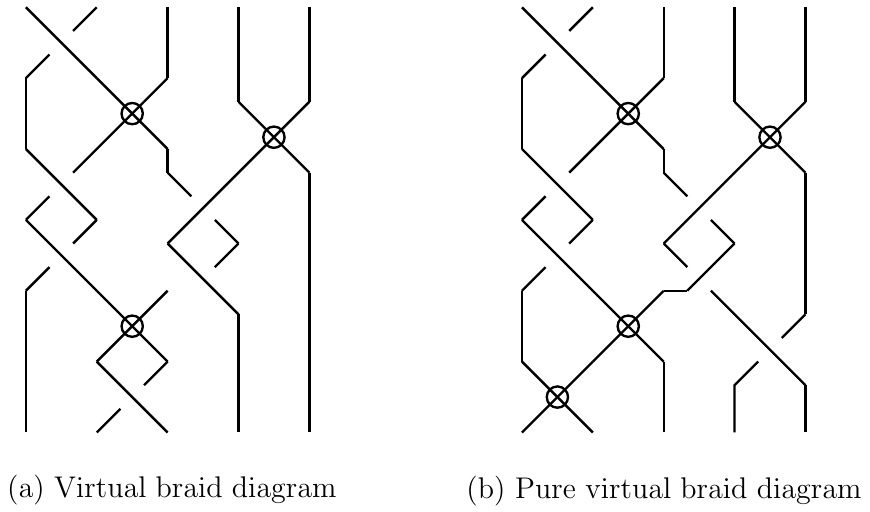}
\caption{Braid diagrams of degree $5$.}
\label{exa}
\end{figure}

Let $n$ be a positive integer. A {\it virtual braid diagram} of degree $n$ is a union of $n$ strands attached to two parallel horizontal bars in ${\mathbb R}^2$, where each strand runs from one bar to the other bar monotonically. Each pair of strands of a virtual braid diagram is allowed to cross in one of three types $\sigma_i, {\sigma_i}^{-1},$ and $v_i$ for some $i\in \{1,\dots,n-1\}$ as depicted in Figure~\ref{fig-sigma} which is called a positive crossing, a negative crossing, and a virtual crossing, respectively. Examples of virtual braid diagrams of degree $5$ are shown in Figure~\ref{exa}.
  We identify two virtual braid diagrams if they are moved to each other by a planar isotopy which fixes the two parallel horizontal bars pointwise.  A {\it classical braid diagram} is a virtual braid diagram which has no virtual crossings.

\begin{figure}[ht]
  \centering
    \includegraphics[width=9cm]{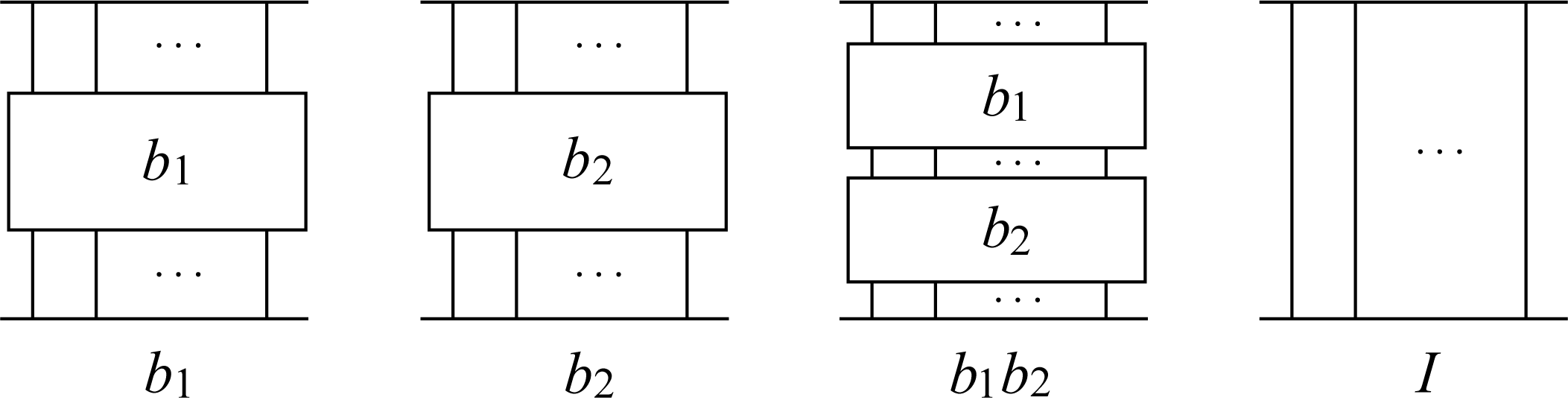}
        \caption{Binary operation for the monoid $BD_n$ or $VBD_n$.}
        \label{con}
        \end{figure}

Let $VBD_n$ (resp. $BD_n$) be the set of all virtual (resp. classical) braid diagrams of degree $n$. For elements $b_1$ and $b_2$ of $VBD_n$, we define an element $b_1b_2$ of $VBD_n$ by placing $b_1$ above $b_2$ so that the bottoms of $b_1$ coincide with the tops of $b_2$ as shown in Figure~\ref{con}. Then, $VBD_n$ forms a monoid by the above binary operation, and is generated by $\sigma_i, {\sigma_i}^{-1}$ and $v_i$ for $i\in \{1,\dots, n-1\}$.  Here, a monoid is a set equipped with an associative binary operation and an identity element.  
The identity element of $VBD_n$ is the diagram without crossings. For example, the virtual braid diagram in Figure~\ref{exa}(a) can be represented by ${\sigma_1}^{-1}v_2v_4{\sigma_1}^{-1}\sigma_3{\sigma_1}^{-1}{\sigma_3}^{-1}v_2{\sigma_2}^{-1}$.  We remark that $BD_n$ is a submonoid of $VBD_n$ generated by $\sigma_i$ and ${\sigma_i}^{-1}$ for $i\in \{1,\dots ,n-1\}$. 
In this paper, a braid diagram means a virtual or a classical braid diagram. 

Each braid diagram $\beta$ has its permutation $\pi$. If a strand is in the $i^{th}$ position at the top and in the $j^{th}$ position at the bottom, we denote $\pi(i)=j$.

\begin{definition}
A pure braid diagram is a braid diagram whose permutation is identity. An example is shown in Figure \ref{exa}(b).
\end{definition}

\begin{definition}
A reducible braid diagram is a braid diagram which has a non-alternating bigon, which is represented by $\sigma_i \sigma_i^{-1}$ or $\sigma_i^{-1}\sigma_i$ for some $i$. 
An irreducible braid diagram is a braid diagram which has no non-alternating bigon. For example, see Figure~\ref{rbd}.
 \begin{figure}[ht]
	\centering
\includegraphics[width=6.5cm]{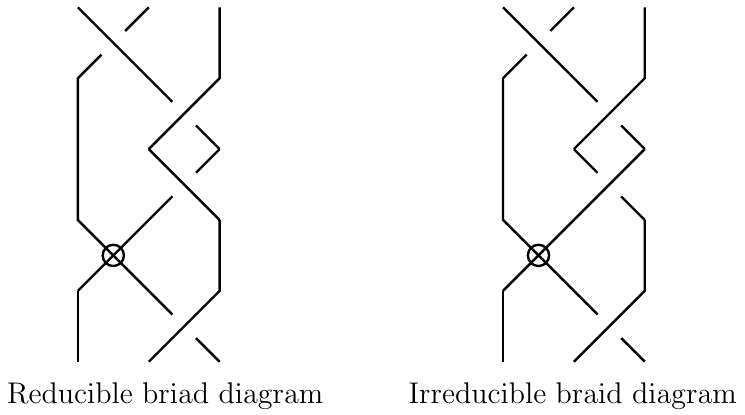}	
  \caption{Reducible and irreducible braid diagrams.}
 \label{rbd}
\end{figure}
\end{definition}

\subsection{Up-down action of braid diagrams}
        
\begin{definition}
    Let $M=(M, \circ )$ be a monoid (a set equipped with an associative binary operation $\circ$ and an identity element \( e \)), and let $\Gamma$ be a set. A right monoid action of \( M \) on \( \Gamma \) is a map
\[\phi: \Gamma \times M \to \Gamma , \quad \text{denoted by } \quad (x, m) \mapsto x \cdot m,\]
such that the following two conditions hold:
\begin{itemize}
    \item For all \( x \in \Gamma \),
   \( x \cdot e = x,\)
   where \( e \) is the identity element of \( M \).

\item For all \( m_1, m_2 \in M \) and \( x \in \Gamma \),
   \( x \cdot (m_1 \circ m_2)  =(x \cdot m_1 )\cdot m_2 . \)
   \end{itemize}
\end{definition}

\begin{definition}
   Let $M$ be a monoid acting on a set $\Gamma$. The orbit of an element $x$ of $\Gamma$ is the set $\{x \cdot m \in \Gamma~|~m\in M\}$, denoted as $x \cdot M$. For any subset $N$ of $M$, we also denote the set $\{ x \cdot m \ | \ m \in N \}$ by $x \cdot N$. 
 \end{definition}

For a braid diagram, each segment of a strand divided by classical crossings is referred to as an {\it edge}. The up-down coloring assigns integers to edges of a braid diagram as follows:  

\begin{definition}[\cite{KAM, KAY}]
The up-down coloring for a braid diagram $\beta$ is a coloring assigning integers to edges of $\beta$ such that each edge which has an over-crossing (resp. under-crossing) at the top point has the integer which is greater (resp. smaller) than the integer of the above edge by one in each strand, as depicted in Figure \ref{pbp2}.
\begin{figure}[ht]
	\centering
\includegraphics[width=4cm]{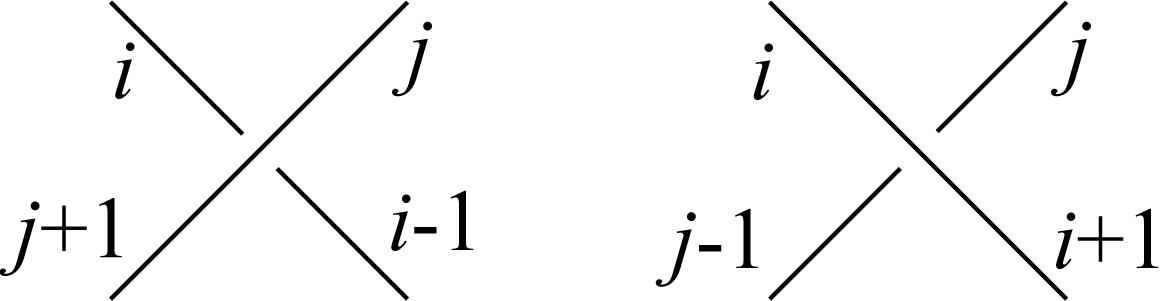}	
  \caption{Up-down coloring.}
 \label{pbp2}
\end{figure}
\label{def-udl}
\end{definition}

\noindent Following the rule, the up-down coloring to a braid diagram $\beta$ of degree $n$ is given uniquely for any initial tuple that consists of integers put to the edges on the top of $\beta$ from left to right. 

\begin{Example}
The up-down coloring for a braid diagram of degree $3$ is shown in Figure \ref{wdll2} for the initial tuples $(0,0,0)$, $(1,2,3)$ and $(x_1, x_2, x_3) \in \mathbb{Z}^3$.
\begin{figure}[ht]
	\centering
\includegraphics[width=6.5cm]{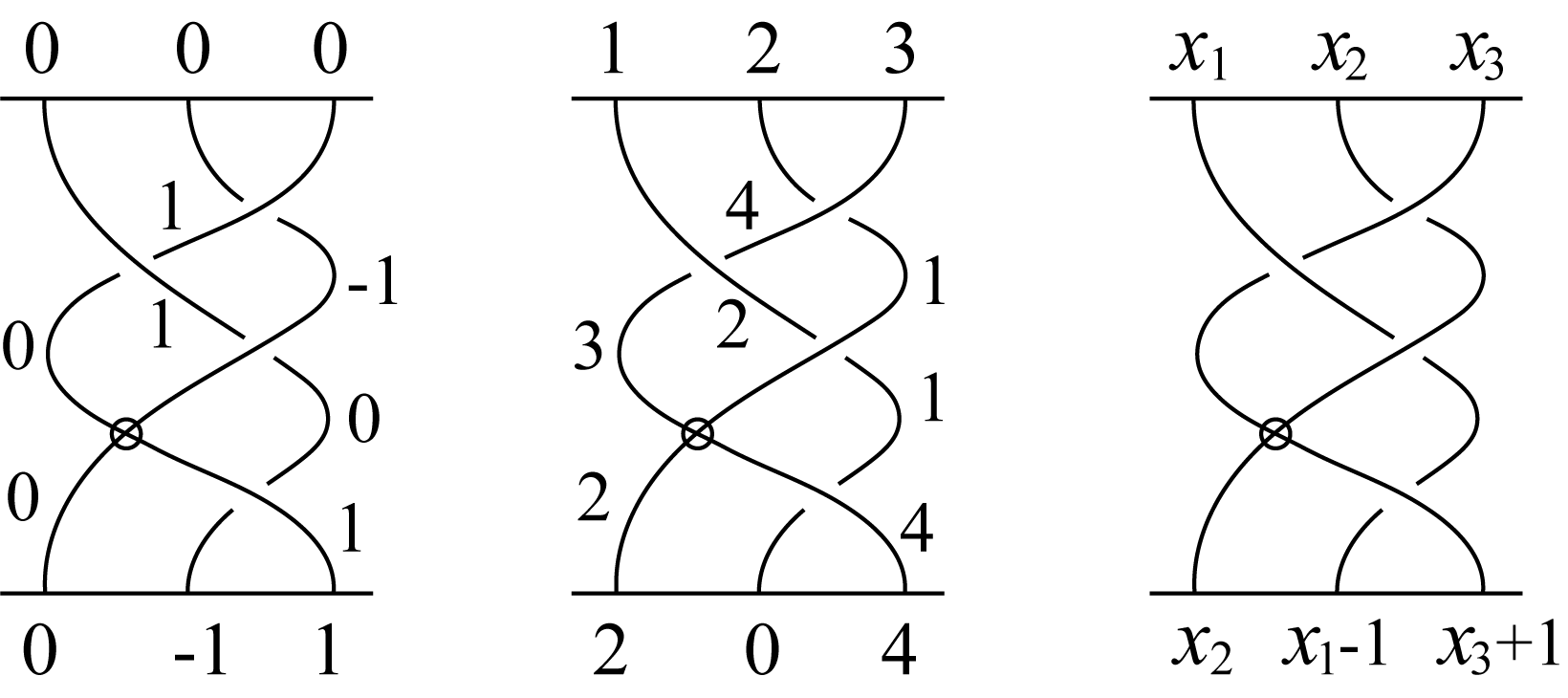}	
  \caption{Up-down colorings for the same braid diagram.}
 \label{wdll2}
\end{figure}
\end{Example}

\begin{figure}[ht]
\centering
\includegraphics[width=2.5cm]{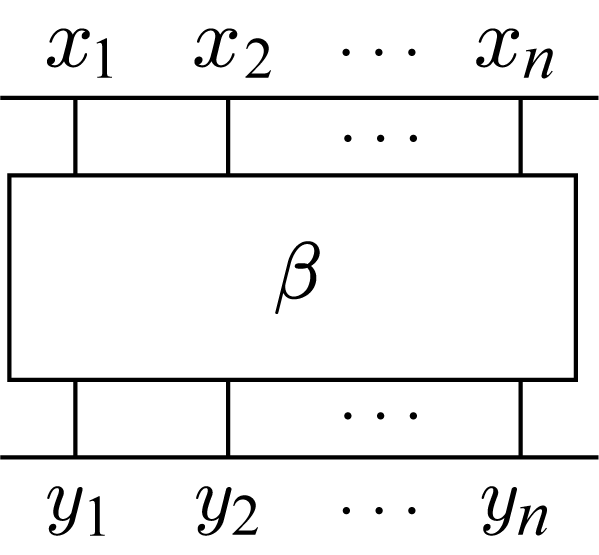}	
  \caption{The up-down action $\vec{x} \cdot \beta = \vec{y}$.}
 \label{func}
\end{figure}

\noindent Let $\beta$ be a braid diagram of degree $n$ with an up-down coloring with initial tuple $\vec{x}=(x_1, x_2, \dots , x_n)$ on the top. When the coloring on the bottom is $\vec{y}=(y_1, y_2, \dots , y_n)$ as shown in Figure \ref{func}, we denote $\vec{x} \cdot \beta = \vec{y}$. 

\begin{prop}\label{monoidVBD}
The map $\mathbb{Z}^n \times VBD_n \to \mathbb{Z}^n$ with $(\vec{x}, \beta ) \mapsto \vec{x} \cdot \beta$ is a right monoid action on $\mathbb{Z}^n$. 
\end{prop}

\begin{corollary}
    The map $\mathbb{Z}^n \times BD_n \to \mathbb{Z}^n$ with $(\vec{x}, \beta ) \mapsto \vec{x} \cdot \beta$ is a right monoid action on $\mathbb{Z}^n$. 
\end{corollary}
\noindent We call the monoid action $\mathbb{Z}^n \times VBD_n \to \mathbb{Z}^n$ with $(\vec{x}, \beta ) \mapsto \vec{x} \cdot \beta$ the {\it up-down action}. 
The up-down action has the geometrical meaning as follows. 
 
\begin{prop}
For each $\beta \in VBD_n$ and $\vec{y}=(y_1, y_2, \dots , y_n) =\vec{0}\cdot \beta$, the value of $y_i$ is the number of over-crossings minus the number of under-crossings on the strand which is placed in the $i^{th}$ position at the bottom. 
\label{prop-o-u}
\end{prop}

\begin{Example}
For the diagram $\beta$ depicted in Figure \ref{wdll2}, we have $(0,0,0) \cdot \beta =(0,-1, 1)$, $(1,2,3) \cdot \beta =(2,0,4)$, and in general $(x_1, x_2, x_3) \cdot \beta =(x_2, x_1-1, x_3+1)$. 
The strand which is in the 3rd position at the bottom has two over-crossings and one under-crossings. Namely, the difference is one. We can also see this from the 3rd component of $\vec{0}\cdot \beta$. 
\label{ex-tf}
\end{Example}

\begin{definition}
A pair of two elements $\vec{x}, \vec{y} \in \mathbb{Z}^n$ are $VBD_n$-equivalent if there exists a virtual braid diagram $\beta \in VBD_n$ such that $\vec{x}\cdot \beta =\vec{y}$, and denoted by $\vec{x} \stackrel{{ \mathrm{VBD_n}}}{\backsim} \vec{y}$. Two elements $\vec{x}, \vec{y} \in \mathbb{Z}^n$ are $BD_n$-equivalent if there exists a classical braid diagram $\beta \in BD_n$ such that $\vec{x}\cdot \beta =\vec{y}$, and denoted by $\vec{x} \stackrel{{ \mathrm{BD_n}}}{\backsim} \vec{y}$.
\end{definition}

\begin{prop}
The relation $\stackrel{{ \mathrm{VBD_n}}}{\backsim}$ is an equivalence relation.
\end{prop}
\begin{proof}
\begin{itemize}
\item[(i)] Reflexive Property: $\vec{x} \stackrel{{ \mathrm{VBD_n}}}{\backsim} \vec{x}$ holds for any $\vec{x}$ by taking the trivial braid diagram $\beta$. 
\item[(ii)] Symmetric Property: For $\vec{x}, \vec{y} \in \mathbb{Z}^n$, if there exists a virtual braid diagram $\beta$ such that $\vec{x}\cdot \beta =\vec{y}$, then there also exists a braid diagram $\beta'$ such that $\vec{y}\cdot \beta' =\vec{x}$. Let $\beta'$ be the braid diagram obtained from $\beta$ by a horizontal rotation, namely, $\beta'=a_n \dots a_2 a_1$ when $\beta = a_1 a_2 \dots a_n$, where each $a_k$ represents a classical or virtual crossing $\sigma_i$, $\sigma_i^{-1}$ or $v_i$ for some $i$. Since each edge of $\beta'$ has the same integer to the corresponding edge of $\beta$ with the up-down coloring starting with $\vec{x}$, $\vec{y}\cdot \beta' =\vec{x}$ holds. Therefore, $\vec{y} \stackrel{{ \mathrm{VBD_n}}}{\backsim} \vec{x}$ holds when $\vec{x} \stackrel{{ \mathrm{VBD_n}}}{\backsim} \vec{y}$.
\item[(iii)] Transitive Property: It follows from Proposition~\ref{monoidVBD}.
\end{itemize}
\end{proof}

\begin{corollary}
The relation $\stackrel{{ \mathrm{BD_n}}}{\backsim}$ is an equivalence relation.
\end{corollary}

%%%%%%%%%%%%%%%%%
\subsection{Properties of up-down action}

In this subsection, we show properties of up-down action which will be used in this paper. 
Recall that $\mathrm{tr}(\vec{x})=\sum_{i=1}^n x_i$ for $\vec{x}=(x_1, x_2, \dots , x_n) \in \mathbb{Z}^n$.

\begin{lemma}
If $\vec{x} \stackrel{{ \mathrm{VBD_n}}}{\backsim} \vec{y}$, then $\mathrm{tr}(\vec{x}) = \mathrm{tr}(\vec{y})$. 
\label{xy-sum}
\end{lemma}

\begin{proof}
Let $\beta \in VBD_n$ be a braid diagram with an up-down coloring. The horizontal sums of the integers of the $n$ edges are equivalent above and below a classical crossing by Definition \ref{def-udl}; if one passes through an over-arc, the integer increases by one and if one passes through an under-arc, the integer decreases by one at a classical crossing. 
There is no change around a virtual crossing. Hence, the horizontal sum of the integers is always constant. Therefore, $\mathrm{tr}(\vec{x})=\mathrm{tr}(\vec{y})$. 
\end{proof}

\begin{prop}
Let $\beta_1$, $\beta_2$ be virtual braid diagrams with $\vec{0}\cdot \beta_1 =(y_1, y_2, \dots , y_n)$, $\vec{0}\cdot \beta_2 =(z_1, z_2, \dots , z_n)$. When $\beta_2$ has the permutation $\pi$, we have \\
$\vec{0}\cdot (\beta_1 \beta_2)=(y_{\pi^{-1}(1)}+z_1, y_{\pi^{-1}(2)}+z_2 , \dots , y_{\pi^{-1}(n)}+z_n)$.
\label{perm-add}
\end{prop}

\begin{proof}
For each braid diagram $\beta$ with $\vec{0}\cdot \beta =(y_1, y_2, \dots , y_n)$, $y_k$ represents the value of the number of over-crossings minus the number of under-crossings on the strand that is placed at the $k^{th}$ position on the bottom by Proposition \ref{prop-o-u}. For $\beta_1 \beta_2$, the strand in $\beta_2$ that is placed in the $k^{th}$ position at the bottom is connected to the strand in $\beta_1$ that is in the $\pi^{-1}(k)^{th}$ position at the bottom of $\beta_1$. Hence, the value of the total number of over-crossings minus the total number of the under-crossings on the strand is $y_{\pi^{-1}(k)}+z_k$. 
\end{proof}

\begin{corollary}
For $\beta_1 , \beta_2 \in VBD_n$, if $\beta_2$ is a pure braid diagram, then $\vec{0}\cdot (\beta_1 \beta_2) = \vec{0}\cdot \beta_1 + \vec{0}\cdot \beta_2$.
\label{pure-add}
\end{corollary}

\noindent In the same way to Proposition \ref{perm-add}, the following lemma follows from Proposition \ref{prop-o-u}. 

\begin{lemma}
For $\vec{x}=(x_1, x_2, \dots , x_n), \ \vec{y}=(y_1, y_2, \dots , y_n) \in \mathbb{Z}^n$ and a braid diagram $\beta$ with permutation $\pi$, we have $\vec{x}\cdot \beta =\vec{y}$ if and only if $\vec{0}\cdot \beta =(y_1-x_{\pi^{-1}(1)}, y_2-x_{\pi^{-1}(2)}, \dots , y_n-x_{\pi^{-1}(n)})$. 
\label{cor-add0}
\end{lemma}

%%%%%%%%%%%%%%%%%%%%%%%%%%%%%%%%%%%
\section{Classical braid diagrams action on $\mathbb{Z}^n$ by up-down coloring}

In this section, our aim is to find a necessary and sufficient condition on $\vec{y}$ such that $\vec{y}$ belongs to the orbit $\vec{x} \cdot BD_n$ of $\vec{x} \in \mathbb{Z}^n$ under the up-down action of $BD_n$. This means we need to find the images of $\vec{x}$ under all the classical braid diagrams. 

\subsection{Preparation}

In this subsection, we show properties of up-down action of classical braid diagrams which will be used in this section. 

\begin{prop}\label{necessary}
    If $\vec{y} \in \vec{0} \cdot BD_n$, then $\mathrm{tr}(\vec{y})=0$.
\end{prop}
\begin{proof}
It follows from Lemma \ref{xy-sum}.
\end{proof}

\begin{lemma}
For each $\beta \in BD_n$ and $\vec{y}=(y_1, y_2, \dots , y_n)=\vec{0}\cdot \beta$, the parity of $y_k$ is equivalent to the parity of the number of crossings on the $\pi^{-1}(k)^{th}$strand, where $\pi$ is the permutation of $\beta$.
\label{yk-p}
\end{lemma}

\begin{proof}
Let $s$ be the strand of $\beta$ which is in the $k^{th}$ position at the bottom. Let $O_s$ (resp. $U_s$) be the number of over-crossings (resp. under-crossings) on the strand $s$. Then, $y_k= O_s - U_s$ by Proposition \ref{prop-o-u}, and the number of crossings on $s$ is $O_s + U_s$. Hence they have the same parity since $O_s-U_s \equiv (O_s-U_s)+2U_s=O_s+U_s \pmod{2}$. 
\end{proof}

\begin{corollary}~\label{xy-pure}
If $\vec{x}=(x_1, x_2, \dots , x_n)$ and $\vec{y}=(y_1, y_2, \dots , y_n) \in \mathbb{Z}^n$ satisfy $\vec{x}\cdot b =\vec{y}$ for some classical pure braid diagram $b$, then $x_i \equiv y_i \pmod{2}$ for each $1\leq i\leq n$. In particular, if $\vec{y}=\vec{0}\cdot b$, then each $y_i$ is even.
\end{corollary}

To find explicitly the classical braid diagram $\beta$ of degree $n$ that gives the output $\vec{y}$ when $\beta$ act on $\vec{0}$, it is interesting to find such a diagram $\beta$ which is an irreducible braid diagram. 
If we allow reducible braid diagrams, then the problem might become easy. For instance, for $n=2$ we can construct reducible braid diagrams for every $\vec{y} \in \mathbb{Z}^n$ such that $\mathrm{tr}(\vec{y})=0$; we have $\vec{0} \cdot (\sigma_1\sigma_1^{-1})^c=(-2c, 2c), \ \vec{0} \cdot (\sigma_1\sigma_1^{-1})^c=(2c+1, -2c-1), \ \vec{0} \cdot (\sigma_1^{-1}\sigma_1)^c=(2c, -2c)$ and $\vec{0} \cdot (\sigma_1^{-1}\sigma_1)^c\sigma_1^{-1}=(-2c+1, 2c-1)$ for any $c>0$. Hence, when $n=2$, we obtain the orbit $\vec{0}\cdot BD_2 = \{ (a,-a) \ | \ a \in \mathbb{Z} \}$. 

Let $IBD_n$ be the set of all irreducible classical braid diagrams of degree $n$. We remark here that $IBD_n$ is not a monoid with the braid product in $BD_n$. 
When $n=2$, we have the following. 

\begin{lemma}
$\vec{0} \cdot IBD_2=\{(0,0),(1,-1),(-1,1)\}.$
\end{lemma}

\begin{proof}
Since we do not allow the non-alternating bigons $\sigma_1 \sigma_1^{-1}$ or $\sigma_1^{-1} \sigma_1$, any braid diagram of degree $2$ is either $(\sigma_1)^n$ or $(\sigma_1^{-1})^n$ for some $n$. When $n$ is an even number, we have $\vec{0}\cdot (\sigma_1)^n=\vec{0}\cdot (\sigma_1^{-1})^n=(0,0)$. When $n$ is an odd number, we have $\vec{0}\cdot (\sigma_1)^n=(1,-1)$, $\vec{0}\cdot (\sigma_1^{-1})^n=(-1,1)$. 
\end{proof}

%%%%%%%%%%%%%%%%%%%%%%%%%%%%%%%%%%%%%%%%%%%%%%%
\subsection{The orbit of $\vec{0} \cdot IPBD_n$}
Let $PBD_n$ denote the set of all classical pure braid diagrams of degree $n$ and let $IPBD_n$ be the set of all irreducible classical pure braid diagrams of degree $n$. When $n=2$, we have 
$PBD_2 = \{ \sigma_1^{\varepsilon_1}\sigma_1^{\varepsilon_2} \dots \sigma_1^{\varepsilon_k} \ | \ \varepsilon_1, \varepsilon_2, \dots , \varepsilon_k \in \{ \pm 1 \} , \ \sum_{i=1}^k \varepsilon_i \equiv 0 \pmod{2} \}$, 
$IPBD_2 = \{ \sigma_1^k \ | \ k \in \mathbb{Z}, \ k \equiv 0 \pmod{2} \}$, and 
$\vec{0}\cdot PBD_2= \{ (2a, -2a) \ | \ a \in \mathbb{Z} \}$, $\vec{0}\cdot IPBD_2 = \{ (0,0) \}$.

\begin{Example}
    For the weaving braid $w=( \sigma_1 \sigma_2^{-1} \sigma_3 \sigma_4^{-1} \dots \sigma_{n-1}^{-(-1)^n})^n$ or $w^*=( \sigma_1^{-1} \sigma_2 \sigma_3^{-1} \sigma_4 \dots \sigma_{n-1}^{(-1)^n})^n \in IPBD_n$ shown in Figure \ref{weaving}, each strand of $w$ or $w^*$ has $2(n-1)$ crossings in the alternating manner. By Proposition \ref{prop-o-u}, both braid diagrams fix the zero element of $\mathbb{Z}^n$ under the up-down action, namely, $ \vec{0} \cdot w =\vec{0}$, and $  \vec{0} \cdot w^* =\vec{0}$.
    \begin{figure}[ht]
	\centering
\includegraphics[width=8cm]{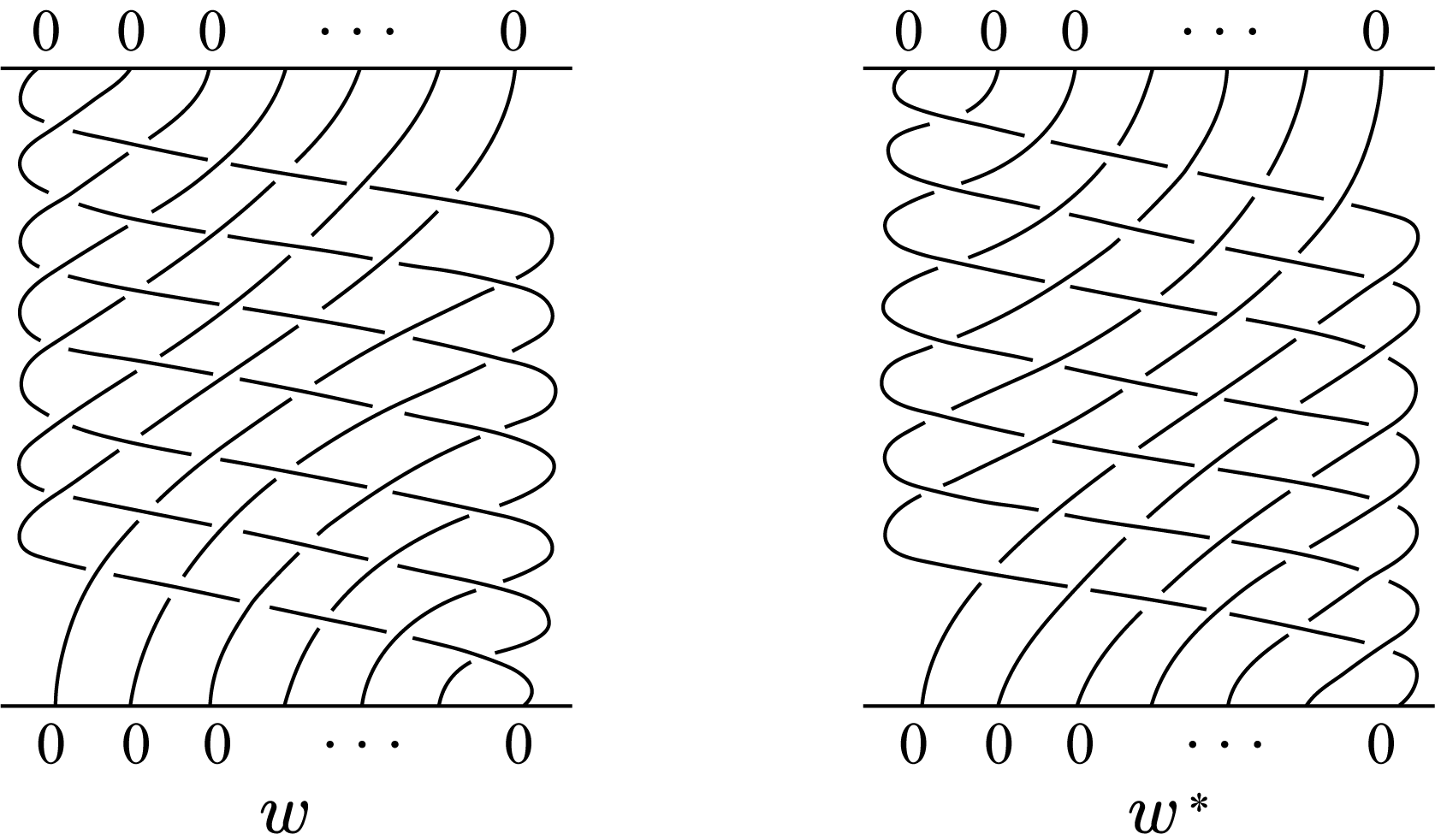}	
  \caption{A weaving braid $w$ and its mirror image $w^*$.}
 \label{weaving}
\end{figure}
\label{ex-weaving}
\end{Example}

\begin{figure}[ht]
	\centering
\includegraphics[width=12cm]{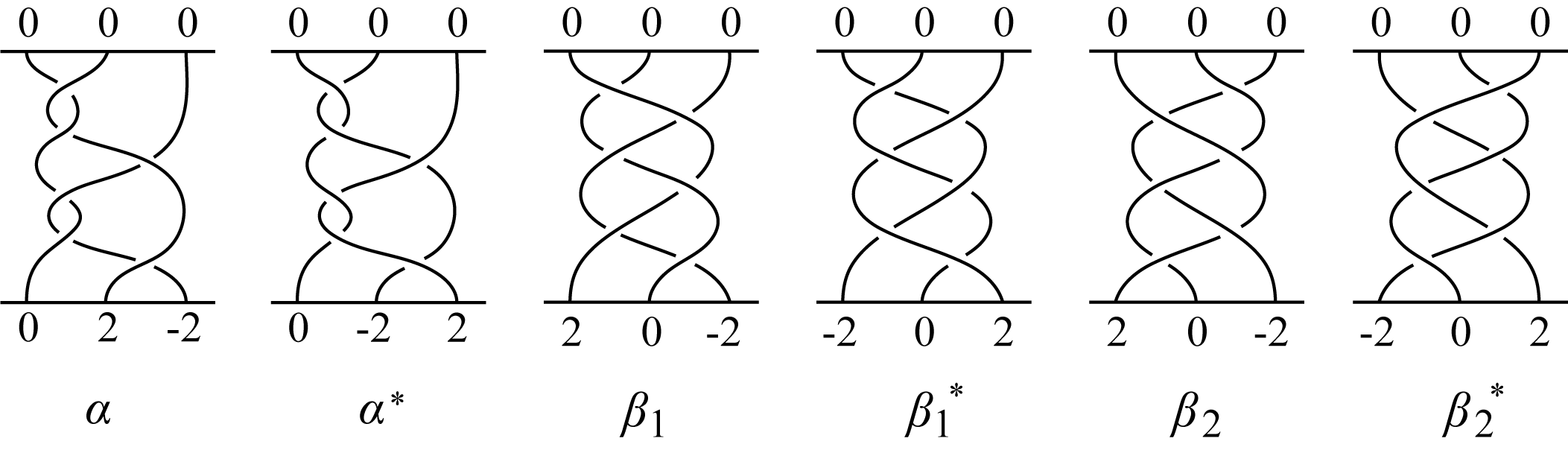}	
  \caption{Irreducible pure braid diagrams of degree 3.}
 \label{Irreducible}
\end{figure}

Let $\alpha, \ \alpha^*, \ \beta_1, \ \beta_1^*, \ \beta_2, \ \beta_2^*$ be the classical pure braid diagrams of degree $3$ as shown in Figure \ref{Irreducible}. Assume that $\alpha^{-k}=( \alpha^* )^k, \ \beta_i^{-k}=( \beta_i^* )^k$ for any $k$ and $i=1, 2$. 

\begin{lemma}
$\vec{0}\cdot IPBD_3 = \{ (2a, 2b, -2a-2b) \ | \ a, b \in \mathbb{Z} \} = \vec{0}\cdot PBD_3$. 
\end{lemma}

\begin{proof}
For each $\beta \in PBD_3$, we have $\vec{0}\cdot \beta =(2a, 2b, -2a-2b)$ for some $a, b \in \mathbb{Z}$ by Proposition \ref{necessary} and Corollary~\ref{xy-pure}. On the other hand, for given $a, b \in \mathbb{Z}$, we can construct a classical braid diagram $\beta$ of degree $3$ such that $\vec{0}\cdot \beta =(2a, 2b, -2a-2b)$ as $\beta= \beta_1^a \alpha^b$ using the classical pure braid diagrams $\alpha, \alpha^*, \beta_1, \beta_1^*$ in Figure \ref{Irreducible} with the rules $\alpha^{-k}=( \alpha^* )^k$ and $\beta_1^{-k}=( \beta_1^* )^k$. The diagram $\beta$ is irreducible for any $a, b \in \mathbb{Z}$. 
\end{proof}

We define an injection $\iota_s^t: BD_n \to BD_{n+s+t}$, where $\iota_s^t(\beta)$ is the classical braid  diagram of degree $n + s + t$ obtained from $\beta \in BD_n$ by adding $s$ trivial strands to the left and $t$ trivial strands to the right.

\begin{theorem}~\label{main1}
    For $\vec{y}=(2a_1,2a_2,\dots,2a_n) \in \mathbb{Z}^n$ ($n \geq 3$) with $2a_n=-2a_1-2a_2 \dots -2a_{n-1}$, there exists an irreducible classical pure braid diagram $b \in IPBD_n$ such that \(  \vec{0} \cdot b = \vec{y} \). Explicitly, the following diagram $b$ satisfies the condition:
    \begin{itemize}
        \item[(i)] If $n$ is even ($n=2k$), then
        \begin{align*}
            b= & \iota_0^{2k-3}(\beta_2^{a_1})\iota_1^{2k-4}(\beta_1^{a_2}) \cdots\cdots \cdots \iota_{2k-5}^{2}(\beta_1^{a_2+a_4\ldots + a_{2k-4}})\iota^1_{2k-4}(\beta_2^{a_1+a_3+\ldots +a_{2k-3}})\\
            & \iota^0_{2k-3}(\beta_1^{a_2+a_4\ldots + a_{2k-2}})\iota^0_{2k-3}(\alpha^{a_1+a_3+\ldots +a_{2k-1}}),
        \end{align*}
        \item[(ii)] If $n$ is odd ($n=2k+1$), then
        \begin{align*}
            b= & \iota_0^{2k-2}(\beta_1^{a_1})\iota_1^{2k-3}(\beta_2^{a_2}) \cdots\cdots \cdots \iota_{2k-4}^{2}(\beta_1^{a_1+a_3\ldots + a_{2k-3}})\iota^1_{2k-3}(\beta_2^{a_2+a_4+\ldots +a_{2k-2}})\\
            & \iota^0_{2k-2}(\beta_1^{a_1+a_3\ldots + a_{2k-1}})\iota^0_{2k-2}(\alpha^{a_2+a_4+\ldots +a_{2k}}),
        \end{align*}
    \end{itemize}
    where $\alpha , \alpha^* , \beta_1 , \beta_1^* , \beta_2 , \beta_2^*$ are the classical pure braid diagrams in Figure \ref{Irreducible}, with $\alpha^{-m}=( \alpha^* )^m, \ \beta_i^{-m}=( \beta_i^* )^m$, for $m>0$.
\end{theorem}

\begin{proof}
We have $\vec{0}\cdot \iota_h^{2k-h-3}(\beta_i^a)=(0,0, \dots ,2a,0,-2a,0, \dots ,0)$ for each $h \in \{ 0, 1,2, \dots ,n-2 \}$, where the $(h+1)^{th}$ component is $2a$. We also have $\vec{0}\cdot \iota_{2k-3}^0(\alpha^c)=(0, \dots , 0,2c,-2c)$. 
Since $\alpha^m$, $(\alpha^*)^m$, $(\beta_i)^m$, $(\beta_i^*)^m$ $(i=1, 2, \ m \in \mathbb{N})$ are pure braids, we have 
\begin{align*}
\vec{0} \cdot b= & \vec{0}\cdot \iota_0^{2k-3}(\beta_2^{a_1})+ \vec{0} \cdot \iota_1^{2k-4}(\beta_1^{a_2}) + \dots \\
& + \vec{0}\cdot \iota^0_{2k-3}(\beta_1^{a_2+a_4\ldots + a_{2k-2}}) + \vec{0}\cdot \iota^0_{2k-3}(\alpha^{a_1+a_3+\ldots +a_{2k-1}})
\end{align*}
when $n$ is even by Corollary \ref{pure-add}. Let $\vec{0}\cdot b=(b_1, b_2, \dots , b_n)$. 
\begin{itemize}
\item The first component $b_1$ of $\vec{0}\cdot b$ is equal to that of $\vec{0}\cdot \iota_0^{2k-3}(\beta_2^{a_1})$ since only $\vec{0}\cdot \iota_0^{2k-3}(\beta_2)$ has a non-zero first component. Hence, $b_1=2a_1$. 
\item The second component $b_2$ of $\vec{0}\cdot b$ is equal to that of $\vec{0}\cdot \iota_1^{2k-4}(\beta_1^{a_2})$ since only $\vec{0}\cdot \iota_1^{2k-4}(\beta_1)$ has a non-zero second component. Hence, $b_2=2a_2$. 
\item When $j$ is odd and $3\leq j \leq n-3=2k-3$, the $j^{th}$ component $b_j$ of $\vec{0}\cdot b$ is that of $\vec{0}\cdot \iota_{j-3}^{2k-j}(\beta_2^{a_1+a_3+ \dots + a_{j-2}})+ \vec{0}\cdot \iota_{j-1}^{2k-2-j}(\beta_2^{a_1+a_3+ \dots + a_{j-2}+a_j})$ since only $\vec{0}\cdot \iota_{j-3}^{2k-j}(\beta_2)$ and $\vec{0}\cdot \iota_{j-1}^{2k-2-j}(\beta_2)$ concern the $j^{th}$ component. Hence, $b_j= -2(a_1+a_3+ \dots +a_{j-2})+2(a_1+a_3+\dots +a_{j-2}+a_j)=2a_j$. 
\item Similarly, when $j$ is even and $4\leq j \leq n-2=2k-2$, the $j^{th}$ component $b_j$ of $\vec{0}\cdot b$ is that of $\vec{0}\cdot \iota_{j-3}^{2k-j}(\beta_1^{a_2+a_4+ \dots + a_{j-2}})+ \vec{0}\cdot \iota_{j-1}^{2k-2-j}(\beta_1^{a_2+a_4+ \dots + a_{j-2}+a_j})$ and $b_j=-2(a_2+a_4+ \dots +a_{j-2})+2(a_2+a_4+\dots +a_{j-2}+a_j)=2a_j$. 
\item The $(n-1)^{th}$ component $b_{n-1}$ of $\vec{0}\cdot b$ is that of $\vec{0}\cdot \iota_{2k-4}^{1}(\beta_2^{a_1+a_3+ \dots + a_{2k-3}})+ \vec{0}\cdot \iota_{2k-3}^{0}(\alpha^{a_1+a_3+ \dots + a_{2k-1}})=-2(a_1+a_3+ \dots +a_{2k-3})+2(a_1+a_3+\dots +a_{2k-3}+a_{2k-1})$, and therefore $b_{n-1}=2a_{2k-1}=2a_{n-1}$. 
\item The $n^{th}$ component $b_n$ is $0-2a_1 -2a_2 - \dots -2a_{n-1}=2a_n$ by Proposition \ref{necessary}. 
\end{itemize}
Therefore, $\vec{0}\cdot b =(2a_1, 2a_2, \dots , 2a_n)$. 
The diagram $b$ is irreducible since there are no possibilities for having $\sigma_i \sigma_i^{-1}$ or $\sigma_i^{-1}\sigma_i$ by the construction of $b$. 
Similarly, it holds for the case that $n$ is odd. 
\end{proof}

\noindent In conclusion, we have 
$$\vec{0}\cdot PBD_n = \{ (2a_1, 2a_2, \dots , 2a_n ) \in \mathbb{Z}^n \ | \ 2a_n = -2a_1-2a_2- \dots -2a_{n-1} \} = \vec{0}\cdot IPBD_n.$$

\subsection{Parity types}

\begin{definition}
For \( \vec{x} = (x_1, x_2, \dots, x_n) \), the component \( x_i \) is classified as Type I if, in the pair \( (i, x_i) \), one of the them is odd while the other is even. The component \( x_i \) is classified as Type II if, in the pair \( (i, x_i) \), both are either odd or even. Let $I(\vec{x})$ (resp. $I\hspace{-1.5pt}I(\vec{x})$) denote the number of components of $\vec{x}$ which are of Type I (resp. Type I\hspace{-0.5pt}I). 
\end{definition}

\begin{Example}
When $\vec{x}=(0,2,1,-1)$, $I(\vec{x})=I\hspace{-1.5pt}I(\vec{x})=2$. When $\vec{y}=(2,1,2,1)$, $I(\vec{y})=4$ and $I\hspace{-1.5pt}I(\vec{y})=0$. When $\vec{z}=(0,1,0,2,1)$, $I(\vec{z})=3$ and $I\hspace{-1.5pt}I(\vec{z})=2$. When $\vec{w}=(1,2,2,2,1)$, $I(\vec{w})=1$ and $I\hspace{-1.5pt}I(\vec{w})=4$. 
\end{Example}

\begin{prop}
If $\vec{y}=\vec{x}\cdot \beta$ for $\beta = \sigma_i^{\pm 1}$, then $I(\vec{x})=I(\vec{y})$.
\label{prop-sigma}
\end{prop}
\begin{proof}
For $\vec{x}=(x_1, x_2, \dots , x_n)$ and $\beta = \sigma_i^{\pm 1}$, we have $\vec{y}=(y_1, y_2, \dots , y_n)=(y_1, \dots , x_{i+1} \pm 1, x_i \mp 1, \dots , y_n)$. This implies that the difference between $\vec{x}$ and $\vec{y}$ are only the $i^{th}$ and $(i+1)^{th}$ components, $y_i =x_{i+1} \pm 1$ and $y_{i+1}=x_i \mp 1$. Since the pairs $y_i$ and $x_{i+1}$; $y_{i+1}$ and $x_i$ have the same types, $\vec{x}$ and $\vec{y}$ have the same number of components of type I and type I\hspace{-0.5pt}I. 
\end{proof}
\noindent Since any classical braid diagram $\beta$ can be represented by $\sigma_i^{\pm 1}$s, we have the following corollary. 
\begin{corollary}
If $\vec{x} \stackrel{{ \mathrm{BD_n}}}{\backsim} \vec{y}$, then $I(\vec{x})=I(\vec{y})$.\label{pres}
\end{corollary}
\begin{remark}
    $I(\vec{x})=I(\vec{y})$ if and only if $I\hspace{-1.5pt}I(\vec{x})=I\hspace{-1.5pt}I(\vec{y})$.
\end{remark}
%%%%%%%%%%%%%%%%%%%%%%
\subsection{The orbit $\vec{x}\cdot IBD_n$ for $\vec{x} \in \mathbb{Z}^n$}

In this section, we prove Theorem \ref{thm-main1}. \\

\noindent \textbf{Theorem~\ref{thm-main1}.} For $\vec{x}=(x_1, x_2, \dots , x_n) \in \mathbb{Z}^n$ ($n \geq 3$), 
$$\vec{x}\cdot BD_n = \left\{ \vec{y} \in \mathbb{Z}^n \ | \ \mathrm{tr}(\vec{x})=\mathrm{tr}(\vec{y}), \ I(\vec{x})=I(\vec{y}) \right\} .$$
Moreover, $\vec{x}\cdot IBD_n=\vec{x}\cdot BD_n$.

\begin{lemma}
For $\vec{x}, \vec{y} \in \mathbb{Z}^n$, $\vec{x} \stackrel{{ \mathrm{BD_n}}}{\backsim} \vec{y}$ if and only if $\mathrm{tr}(\vec{x}) = \mathrm{tr}(\vec{y})$ and $I(\vec{x})=I(\vec{y})$. 
\label{xy-iff}
\end{lemma}

\begin{proof}
($\Rightarrow$) By Lemma \ref{xy-sum} and Corollary \ref{pres}. \\
($\Leftarrow$) When $I(S_x)=I(S_y)$, take any permutation $\pi$ such that $\pi(i)=j$ if $x_i$ and $y_j$ have the same type. Take any classical braid diagram $\gamma$ that has the permutation $\pi$. Let $\vec{z}=\vec{x}\cdot \gamma$. Then $z_i \equiv y_i \pmod{2}$ for any $i$ by the permutation $\pi$ and the proof of Proposition \ref{prop-sigma}. Let $\vec{w}=(w_1, w_2, \dots , w_n)=\vec{y}-\vec{z}$. Then, $w_i$ is an even number for any $i$. Also, $\mathrm{tr}(\vec{w}) =\mathrm{tr}(\vec{y}) - \mathrm{tr}(\vec{z}) =\mathrm{tr}(\vec{x}) - \mathrm{tr}(\vec{z}) =0$. By Theorem \ref{main1}, there exists a classical pure braid $b$ such that $\vec{0}\cdot b=\vec{y}-\vec{z}$. Using Lemma \ref{cor-add0} with $\pi =id$, we also have $\vec{z} \cdot b=\vec{y}$. Hence, by taking $\beta =\gamma b$, we have a classical braid diagram $\beta$ such that $\vec{x} \cdot \beta =\vec{y}$. 
\end{proof}

A {\it permutation braid diagram} is a braid diagram $\beta$ which has no negative crossings and each pair of strands of $\beta$ has at most one crossing. For each permutation $\pi$ of degree $n$, we can take a permutation braid diagram that has $\pi$. 

\begin{lemma}
If $\vec{x} \stackrel{{ \mathrm{BD_n}}}{\backsim} \vec{y}$, then there exists an irreducible classical braid diagram $\beta$ of degree $n$ such that $\vec{x}\cdot \beta =\vec{y}$. 
\label{xy-irreducible}
\end{lemma}

\begin{proof}
For the permutation $\pi$ constructed in the proof of Theorem \ref{xy-iff}, take a permutation braid diagram $\gamma$ and let $\vec{z}=\vec{x}\cdot \gamma$. Take a classical pure braid diagram $b$ such that $\vec{z}\cdot b = \vec{y}$ as that constructed in the proof of Theorem \ref{main1}. Then, we have a required classical braid diagram $\beta$ as follows: 
\begin{itemize}
\item[] Step 1. Take $\beta= \gamma b$. Then $\beta$ satisfies $\vec{x}\cdot \beta=\vec{y}$ by Theorem \ref{xy-iff}. If $\beta$ is already irreducible, we fix $\beta$. If $\beta$ has a non-alternating bigon, then the bigon must be found between $\gamma$ and $b$ since $\gamma$ and $b$ have no such bigons. In this case, we apply the next step. 
\item[] Step 2. When there is a non-alternating bigon between $\gamma$ and $b$, then add a pure weaving braid $w$ between them. It will not affect the final outcome $\vec{y}$ as mentioned in Example \ref{ex-weaving}. We also note that there are no non-alternating bigons between $\gamma$ and $w$ since the crossings of $\gamma$ are all positive crossings and there is only one possibility for the existence of a bigon between them at the top-left of $w$ with the positive crossing. If $\gamma w b$ is irreducible, we fix $\beta = \gamma w b$.
\item[] Step 3. If there is a non-alternating bigon between $w$ and $b$, retake $\beta= \gamma ww^* b$. We note that at least one of $\gamma wb$ and $\gamma ww^*b$ has no non-alternating bigons between $w$ or $w^*$ and $b$ since there is only one possibility for the existence of a bigon between them at the lower-right of $w$ or $w^*$. 
\end{itemize}
\begin{figure}[ht]
\centering
\includegraphics[width=8cm]{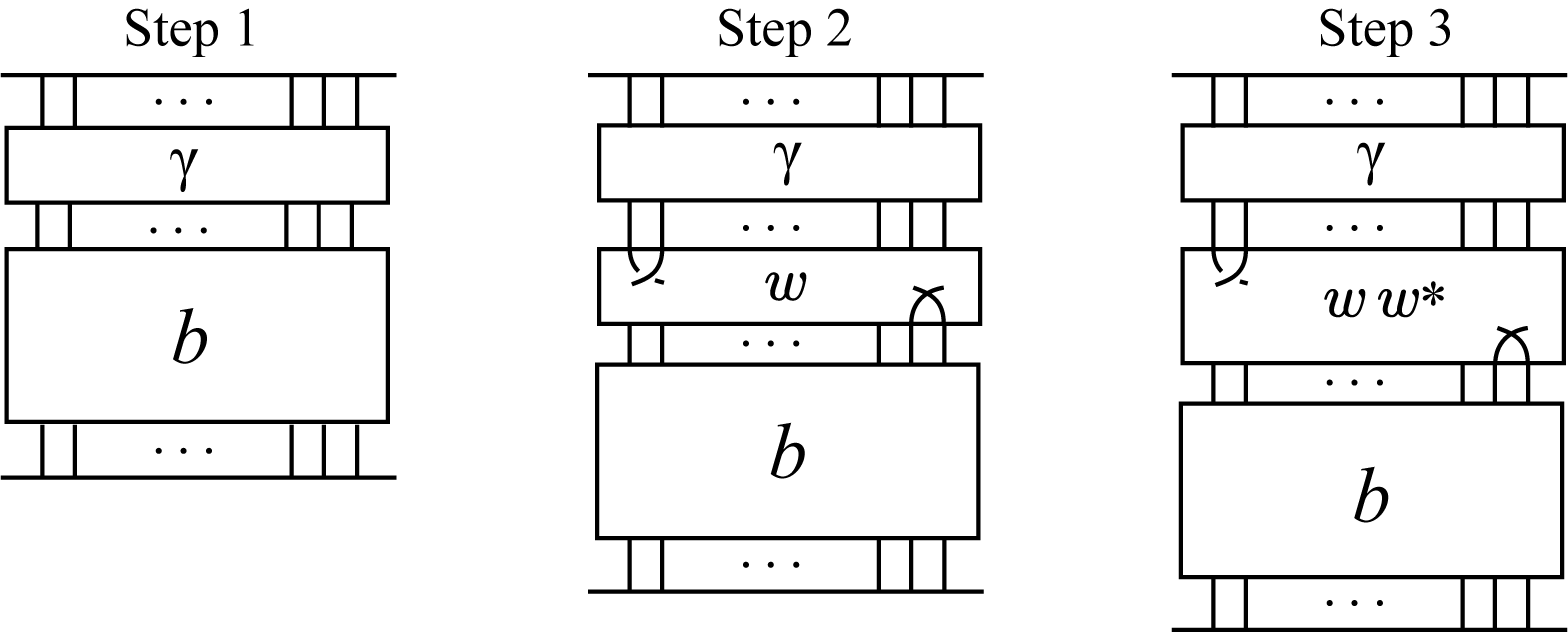}	
\caption{At least one of $\beta =\gamma b, \ \gamma wb, \  \gamma ww^*b$ is irreducible.}
\label{cases}
\end{figure}

\end{proof}
\begin{Example}
Let $\vec{x}=(1,4,-7,-1, 2,9)$, $\vec{y}=(4,0,8,-5,3,-2)$ with $\mathrm{tr}(\vec{x})=\mathrm{tr}(\vec{y})=8$. Since $I(S_x)=I(S_y)=3$, $\vec{x} \stackrel{{ \mathrm{BD_n}}}{\backsim} \vec{y}$. 
Take a permutation $\pi$ and $\gamma$ as shown in Figure \ref{fig-gamma}. 
\begin{figure}[ht]
\centering
\includegraphics[width=2.5cm]{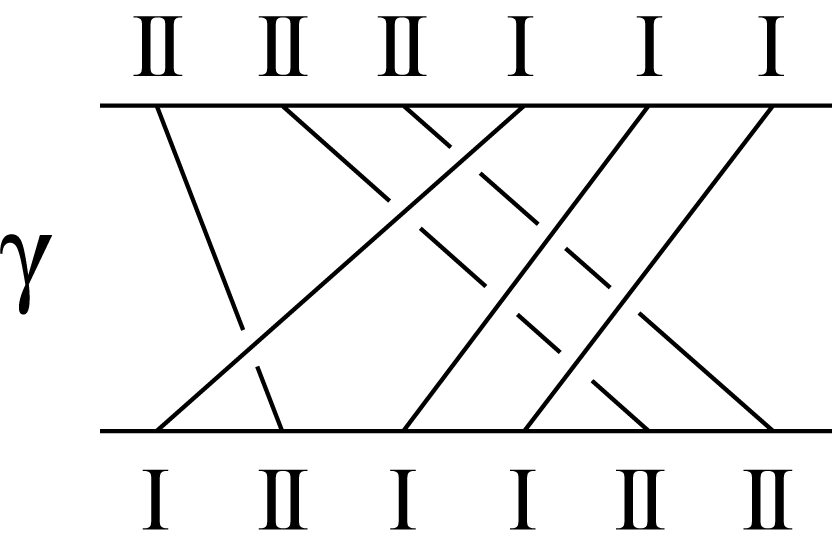}	
\caption{Classical braid diagram $\gamma$ with a permutation $\pi$.}
\label{fig-gamma}
\end{figure}
Then we have $\vec{z}=\vec{x}\cdot \gamma =(2,0,4,11,1,-10)$ and $\vec{y}-\vec{z}=(2,0,4,-16,2,8)$. By the formula of Theorem \ref{main1}, we have an irreducible classical pure braid diagram $b=\iota_0^3 (\beta_2) \iota_2^1(\beta_2^3) \iota_3^0((\beta_1^*)^8) \iota_3^0(\alpha^4)$ satisfying $\vec{0}\cdot b= \vec{y}-\vec{z}$ and equivalently $\vec{z}\cdot b=\vec{y}$. 
Since the diagram $\beta = \gamma b$ has a non-alternating bigon between $\gamma$ and $b$ as shown in Figure \ref{fig-step2}, retake $\beta$ as $\gamma wb$. Thus, we obtain an irreducible diagram $\beta$ such that $\vec{x}\cdot \beta = \vec{y}$. 
\begin{figure}[ht]
\centering
\includegraphics[width=6cm]{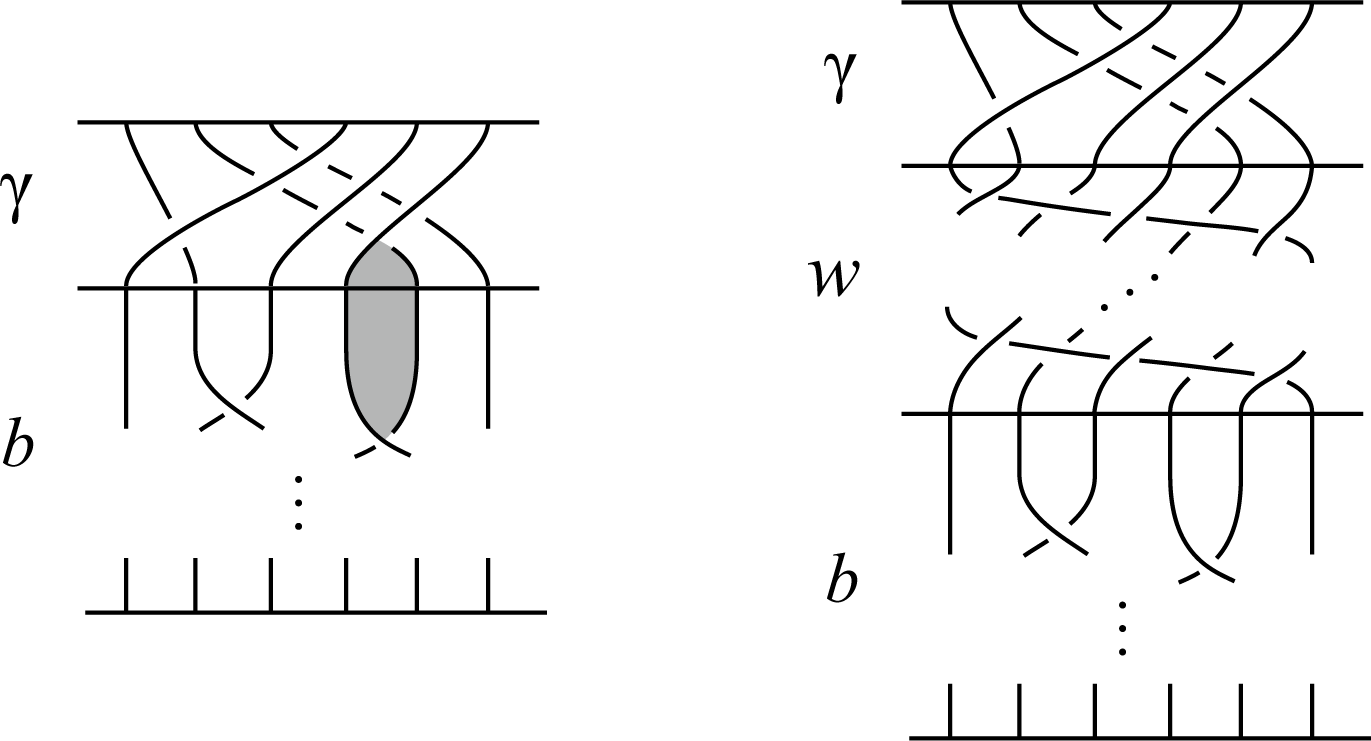}	
\caption{$\gamma b$ is reducible and $\gamma wb$ is irreducible.}
\label{fig-step2}
\end{figure}
\end{Example}

\noindent Now we prove Theorem \ref{thm-main1}. \\

\noindent {\it Proof of Theorem \ref{thm-main1}.} \ 
It follows directly from Lemmas~\ref{xy-iff} and \ref{xy-irreducible}. \qed

\begin{remark}
When $n=2$, 
$$\vec{x}\cdot IBD_2 = \{ (x_1, x_2), (x_2+1, x_1-1), (x_2-1, x_1+1) \} \neq \vec{x}\cdot BD_2$$
for each $\vec{x}=(x_1, x_2)$.
\end{remark}

\begin{corollary}
When $n \geq 3$, we have 
\begin{align*}
\vec{0} \cdot IBD_n = \left\{ \vec{y} \in \mathbb{Z}^n \ | \ \mathrm{tr}(\vec{y})=0 \text{ and } I(\vec{y})= \left\lceil \frac{n}{2} \right\rceil \right\} =\vec{0}\cdot BD_n.
\end{align*}
\end{corollary}

%%%%%%%%%%%%%%%%%%%%%%%%%%%%%%%%%%%%%%%%%%%%%%%%%%%%%%
\section{Virtual braid diagrams action on $\mathbb{Z}^n$ by up-down coloring}

In this section we prove Theorem \ref{thm-main2}. \\

\noindent \textbf{Theorem~\ref{thm-main2}.} For $\vec{x}=(x_1, x_2, \dots , x_n) \in \mathbb{Z}^n$ ($n \geq 3$), 
$$\vec{x}\cdot VBD_n = \left\{ \vec{y} \in \mathbb{Z}^n \ | \ \mathrm{tr}(\vec{x})=\mathrm{tr}(\vec{y}) \right\}.$$ Moreover, $\vec{x}\cdot IVBD_n=\vec{x}\cdot VBD_n$. \\

First, we will determine $\vec{0} \cdot VBD_n$. 
Let $PVBD_n$ be the set of all pure virtual braid diagrams of degree $n$. Let $IPVBD_n$ be a set of all irreducible pure virtual braid diagrams of degree $n$. Then, $IPVBD_n\subseteq PVBD_n$. For $n=2$, we have the following lemma. 
\begin{lemma}~\label{n=2v}
$\vec{0} \cdot PVBD_2=\{ (m,-m)~|~m\in \mathbb{Z}\}.$  Moreover, $\vec{0}\cdot IPVBD_2=\vec{0}\cdot PVBD_2$.
\end{lemma}
\begin{proof}
We have $\vec{0} \cdot PVBD_2 \subseteq \{ (m,-m)~|~m\in \mathbb{Z}\}$ by Lemma \ref{xy-sum}. 
When $m<0$ (resp. $m>0$), we have $\vec{0} \cdot \alpha^m =(m, -m)$ (resp. $\vec{0} \cdot \beta^m =(m,-m)$) for $\alpha^m$, $\beta^m \in IPVBD_2$, where $\alpha$ and $\beta$ are the virtual pure braid diagrams depicted in Figure \ref{vbd}. When $m=0$, take the trivial braid diagram. Hence, $\vec{0} \cdot IPVBD_2 \supseteq \{ (m,-m) \ | \ m \in \mathbb{Z} \}$.

\begin{figure}[ht]
\centering
\includegraphics[width=3cm]{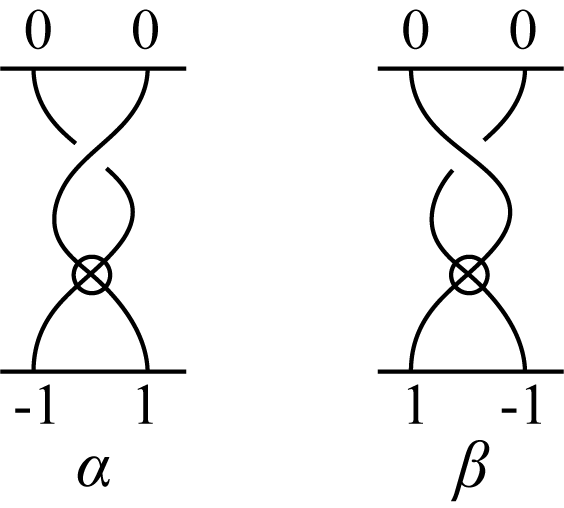}	
  \caption{Virtual pure braid diagrams.}
 \label{vbd}
\end{figure}
\end{proof}
\begin{corollary}
$\vec{0} \cdot VBD_2=\{ (m,-m)~|~m\in \mathbb{Z}\}$.
\end{corollary}
For non-negative integers $s$ and $t$, let $\iota_s^t(\beta)$ be the virtual braid diagram of degree $n + s + t$ obtained from $\beta \in VBD_n$ by adding $s$ trivial strands to the left and $t$ trivial strands to the right. This induces an injection $\iota_s^t: VBD_n \to VBD_{n+s+t}$. 
\begin{lemma}~\label{main3}
    For $\vec{y}=(a_1,a_2,\dots,a_n) \in \mathbb{Z}^n$ with $a_n=-a_1-a_2 \dots -a_{n-1}$, there exists an irreducible pure virtual braid diagram $b \in IPVBD_n$ such that \(  \vec{0} \cdot b = \vec{y} \). Explicitly, the following braid diagram $b$ satisfies the condition:
    $$b=  \iota_0^{n-2}(\gamma_1^{c_1})\iota_1^{n-3}(\gamma_2^{c_2}) \dots \iota_{n-2}^{0}(\gamma_{n-1}^{c_{n-1}}),$$
 where $c_i=~|a_1+a_2+ \dots + a_i|~$, $\gamma_i=\alpha$ when $a_1 +a_2 + \dots +a_i <0$, $\gamma_i=\beta$ when $a_1 +a_2 + \dots +a_i >0$, and $\gamma_i=I$ when $a_1 +a_2 + \dots +a_i =0$ for $1\leq i\leq n-1$.
\end{lemma}
\begin{proof} Let $c>0$. We have $$\vec{0}\cdot \iota_h^{n-h-2}(\beta^c)=(0,0, \dots ,0,c,-c,0, \dots ,0)$$ and $$\vec{0}\cdot \iota_h^{n-h-2}(\alpha^c)=(0,0, \dots ,0,-c,c,0, \dots ,0)$$ for each $h \in \{ 1,2, \dots ,n-2 \}$, where the $(h+1)^{th}$ component is the first non-zero component from the left.  
Observe that the construction of $b$ is irreducible virtual braid, since $b$ contains classical and virtual crossing in the alternating manner.
\end{proof}

\begin{corollary}
$\vec{0} \cdot IVBD_n=\{ \vec{y} \in \mathbb{Z}^n~|~\mathrm{tr}(\vec{y})=0\} = \vec{0}\cdot VBD_n$.
\end{corollary}
\begin{proof}
It follows from Lemmas~\ref{xy-sum} and \ref{main3}.
\end{proof}

\noindent We prove Theorem \ref{thm-main2}.\\

\noindent {\it Proof of Theorem \ref{thm-main2}.} \ 
($\subseteq$) By Lemma~\ref{xy-sum}. \\
($\supseteq$)
Consider $\vec{z}=\vec{y}-\vec{x}$. Since $\mathrm{tr}(\vec{z})=\mathrm{tr}(\vec{y})-\mathrm{tr}(\vec{x})=0$, there exists a virtual braid diagram $b \in IPVBD_n\subseteq IVBD_n\subseteq VBD_n$ such that \(  \vec{0} \cdot b = \vec{z} \) by Lemma~\ref{main3}. 
Since $b$ is a pure virtual braid diagram, $\vec{x} \cdot b =\vec{z}+\vec{x}=\vec{y}$.
\qed

\begin{remark}
For each $\vec{x} \in \mathbb{Z}^n$, $\vec{x} \cdot IPVBD_n=\vec{x} \cdot IVBD_n=\vec{x} \cdot PVBD_n=\vec{x} \cdot VBD_n$.
\end{remark}

\begin{Example}
Consider $\vec{x}=(8,4,2,2)$ and $\vec{y}=(2,3,5,6)$, which satisfy $\mathrm{tr}(\vec{y})=\mathrm{tr}(\vec{x})$. We can construct an irreducible virtual braid diagram $b$ such that \(  \vec{x} \cdot b = \vec{y} \).
Consider $\vec{z}=\vec{y}-\vec{x}=(-6,-1,3,4)$. 
Since $\mathrm{tr}(\vec{z})=0$, we have the virtual braid diagram $b=\iota_0^{2}(\alpha^{6})\iota_1^{1}(\alpha^{7})\iota_2^{0}(\alpha^{4})$ such that $\vec{0} \cdot b = \vec{z}$ by Lemma~\ref{main3}. 
Therefore, \(\vec{x} \cdot b = \vec{y} \).
\end{Example}

\section{The OU matrix and isotropy submonoid by up-down coloring}

The OU matrix, defined in \cite{AY}, represents the relative over/under information for each pair of strands of a classical braid diagram. In this paper, we also define the OU matrix for virtual braid diagrams. 

\begin{definition}
The OU matrix $M(b)$ of $b \in VBD_n$ is an $n\times n$ matrix such that the components on the main diagonal are all zero and each $(i,j)$-component has the number of crossings on the $i^{th}$ strand which are over the $j^{th}$ strand, where $k^{th}$ strand means the strand placed in the $k^{th}$ position at the top. 
\label{def-oum1}
\end{definition}

\begin{Example}
The classical braid diagram $b$ in Figure \ref{ex-oum} has the OU matrix $M(b)$ as the table on the right-hand side. 
\begin{figure}[ht]
	\centering
\includegraphics[width=5cm]{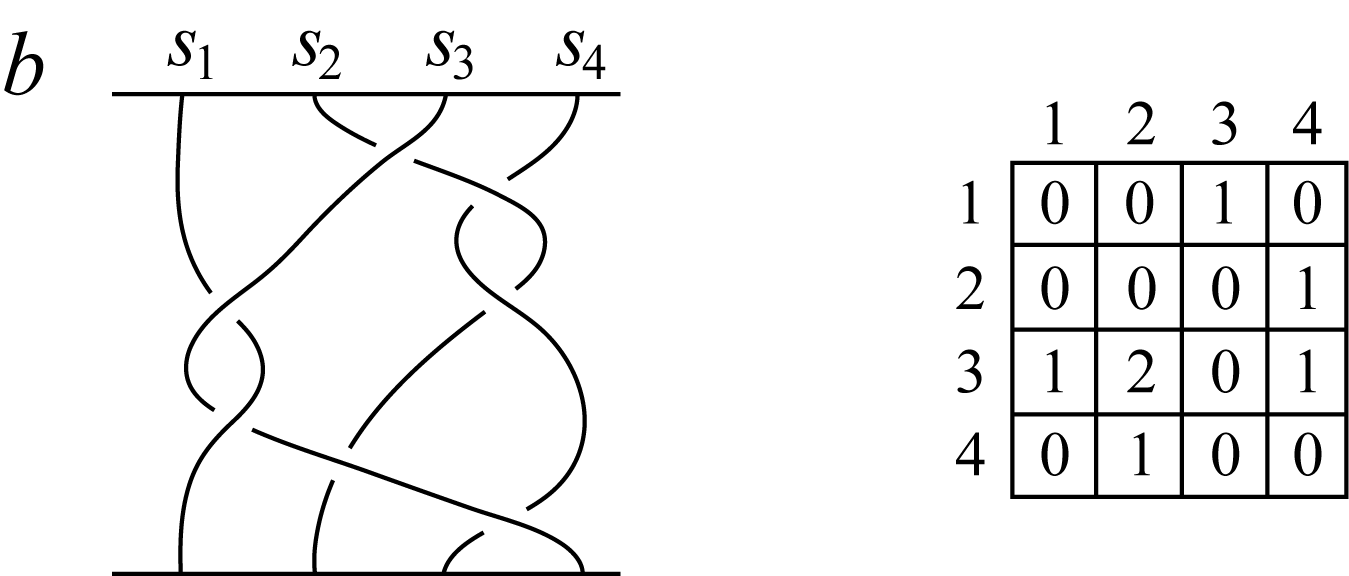}	
  \caption{A classical braid diagram $b$ and its OU matrix.}
 \label{ex-oum}
\end{figure}

\end{Example}
Suppose the OU matrix $M(b)$ of a virtual braid diagram $b$. Let $M_{ij}(b)$ be the $(i,j)$-component of the matrix $M(b)$.
\begin{prop}~\label{ouf}
   Let $b$ be a virtual braid diagram with the permutation $\pi$. Let $r_i=\sum_{j=1}^{n} M_{ij}(b)$ and $c_j=\sum_{i=1}^{n} M_{ij}(b)$. Then, $\vec{0}\cdot b =(y_1,y_2,\ldots, y_n)$ satisfies $y_{\pi (i)}=r_i-c_i$ for all $1\leq i\leq n$. In general, if $(x_1,x_2,\ldots,x_n) \cdot b=(y_1,y_2,\ldots, y_n)$, then $y_{\pi(i)}=x_i+(r_i-c_i)$ for all $1\leq i\leq n$.
\end{prop}
\begin{proof}
By the definition of the OU matrix, $r_i$ (resp. $c_i$) implies the number of over-crossings (resp. under-crossings) on the $i^{th}$ strand. 
\end{proof}

From Proposition~\ref{ouf}, we can also calculate the up-down action using the OU matrix for virtual braid diagrams.
\begin{Example}
    We can calculate the up-down action of the braid diagram $b$ given in Figure~\ref{ex-oum}. The permutation of $b$ is $\pi =(3~~4)$. Here, $r_1=1,~r_2=1,~r_3=4,~r_4=1$, and $c_1=1,~c_2=3,~c_3=1,~c_4=2$. Hence, 
    \begin{align*}
        (x_1,x_2,x_3,x_4) \cdot \beta & =(x_1+(1-1),~x_2+(1-3),~x_4+(1-2),~x_3+(4-1))\\
        &=(x_1,~x_2-2,~x_4-1,~x_3+3).
    \end{align*}
\end{Example}

Let $I\hspace{-1.5pt}S_n$ be the isotropy submonoid of the monoid action on $VBD_n$, namely the set of all virtual braid diagrams $b$ of degree $n$ such that $\vec{0}\cdot b =\vec{0}$. Then, the following holds. 

\begin{corollary}
A virtual braid diagram $b \in VBD_n$ belongs to $I\hspace{-1.5pt}S_n$ if and only if $\sum_{j=1}^n M_{kj}(b)=\sum_{i=1}^n M_{ik}(b)$ for any $k \in \{ 1,2, \dots ,n \}$ on the OU matrix $M(b)$ of $b$. 
\end{corollary}

\noindent It is shown in \cite{AY} that $M(b)$ is a symmetric matrix for any classical positive pure braid diagram $b$. Hence, any classical positive pure braid diagram of degree $n$ belongs to $I\hspace{-1.5pt}S_n$. Any alternating classical braid diagram such that each strand has an even number of crossings also belongs to $I\hspace{-1.5pt}S_n$. 

The up-down coloring for oriented virtual link diagrams was defined in \cite{KAY} to be the coloring of integers to the edges so that if one passes through an over-crossing (resp. under-crossing), the integer increases (resp. decreases) by one. It was shown in \cite{KAY} that there are link diagrams such that giving up-down coloring is impossible. Since the closure of a braid diagram $\beta$ admits an up-down coloring if and only $\beta$ admits an up-down coloring which has the same integers on the top and bottom, we have the following proposition. 

\begin{prop}
The closure of a virtual braid diagram $b$ of degree $n$ admits an up-down coloring if and only if $b$ satisfies $\vec{x}\cdot b=\vec{x}$ for some $\vec{x} \in \mathbb{Z}^n$. 
\label{prop-link}
\end{prop}

\noindent We have the following corollary from Propositions \ref{ouf} and \ref{prop-link}.

\begin{corollary}
If a virtual braid diagram $b$ satisfies $\sum_{j=1}^n M_{kj}(b_{\mathbf{s}})=\sum_{i=1}^n M_{ik}(b_{\mathbf{s}})$ for any $k \in \{ 1,2, \dots ,n \}$, then the closure of $b$ admits an up-down coloring. 
\end{corollary}
\begin{proof}
In this case, $\vec{x} \cdot b=\vec{x}$ holds for any $\vec{x}=(x,x, \dots ,x)$. 
\end{proof}

\section{Torus and Weaving braids}
With the help of the OU matrix, we calculate the up-down action under two braid classes, torus braids and weaving braids.
\begin{prop}
Let $p,q>0$. Consider the torus braid diagram 
$$ T(p,q)= (\sigma_1 \sigma_2 \dots \sigma_{p-1})^{q}.$$
Let $M(T(p,q))$ be the OU matrix of the torus braid diagram $T(p,q)$. Let $r_i$ (resp. $c_j$) be the sum of all the components in the $i^{th}$ row (resp. $j^{th}$ column) of $M(T(p,q))$. Then, the following condition holds:
    \begin{itemize}
        \item[(i)] If $p\geq q$, then
        
        \begin{minipage}{0.4\textwidth}
$$r_i = \begin{cases}
~q-1  &  (i \leq q) \\
~q  & (i > q)
    \end{cases}$$
\end{minipage}
\begin{minipage}{0.4\textwidth}
$$c_i = \begin{cases}
~p-1  &  (i \leq q) \\
~0  & (i > q)
    \end{cases}$$
    
\end{minipage}

Hence, we can conclude  for $p\geq q$. If $\vec{0}\cdot  T(p,q)=\vec{y}=(y_1,\dots, y_p)$, then $\vec{y}$ satisfies the following:
        $$y_{\pi(i)} = \begin{cases}
~q-p  &  (i \leq q) \\
~q  & (i > q)
    \end{cases}$$
    where $\pi$ is the permutation of the torus braid diagram $T(p,q)$.
        \item[(ii)] If $p< q$, let $q=a_1p+a_2$. Then
        
        \begin{minipage}{0.6\textwidth}
$$r_i = \begin{cases}
~(a_2-1)(a_1+1)+(p-a_2)a_1  &  (i \leq a_2) \\
~a_2(a_1+1)+(p-a_2-1)a_1  & (i > a_2)
    \end{cases}$$
\end{minipage}
\begin{minipage}{0.3\textwidth}
$$c_i = \begin{cases}
~(a_1+1)(p-1)  &  (i \leq a_2) \\
~a_1(p-1)  & (i > a_2)
    \end{cases}$$
\end{minipage}

Hence, we can conclude  for $p< q$. If $\vec{0}\cdot  T(p,q)=\vec{y}=(y_1,\dots, y_p)$, then $\vec{y}$ satisfies the following:
        $$y_{\pi(i)} = \begin{cases}
~a_2-p  &  (i \leq a_2) \\
~a_2  & (i > a_2)
    \end{cases}$$
    where $\pi$ is the permutation of the torus braid $T(p,q)$.
    \end{itemize}
\end{prop}

\begin{corollary}
Let $p,q>0$. Consider the weaving braid
   $$ W(p,q)= (\sigma_1 \sigma_2^{-1} \dots \sigma_{p-1}^{-(-1)^p})^{q}.$$
    If $p\geq q$ and $\vec{0}\cdot W(p,q)=\vec{y}=(y_1,\dots, y_p)$, then $\vec{y}$ satisfies the following:
    \begin{itemize}
        \item[(i)] If $p$ and $q$ are even, then $y_k=0$ for all $1\leq k\leq n$.
        \item[(ii)] If $p$ is even and $q$ is odd, then $y_k=1$ for $k\equiv\pi(1) \pmod{2}$, otherwise $y_k=-1$.
        \item[(iii)] If $p$ is odd and $q$ is even, then $y_k=0$ for all $k< \pi(1)$, and for all $k\geq \pi(1)$, $y_k=1$, if $k\equiv \pi(1) \pmod{2}$, otherwise $y_k=-1$.
        \item[(iv)] If $p$ and $q$ are odd, then $y_k=0$ for all $k\geq \pi(1)$, and for all $k< \pi(1)$, $y_k=-1$, if $k\equiv 1 \pmod{2}$, otherwise $y_k=-1$,
    \end{itemize}
    where $\pi$ is the permutation of the weaving braid $W(p,q)$.
\end{corollary}
%%%%%%%%%%%%%%%%%%%%%%%%%%%%%%%
%%%%%%%%%%%%%%%%%%%%%%%%%%%%%%%%%%%%%%%%%%%
\section*{Conclusion}
In this paper, we primarily focused on exploring the orbits generated by the up-down coloring for braid diagrams. This investigation was carried out in two stages: first, we analyzed the orbit \( \vec{x} \cdot BD_n \), and followed by \( \vec{x} \cdot VBD_n \). Future work may involve extending this study to find the orbit under the action of twisted virtual braid diagrams, that is, \( \vec{x} \cdot TVBD_n \). This problem will be quite challenging and interesting because for twisted virtual braid diagrams the up-down coloring behaves differently.

\section*{Acknowledgements}
The first author would like to thank the University Grants Commission(UGC), India, for Research Fellowship with NTA Ref.No.191620008047. 
The second author's work was partially supported by JSPS KAKENHI Grant Number JP21K03263. 
The third author's work was partially supported by JSPS KAKENHI Grant Number JP19K03508. 
The fourth author acknowledges the support of the Anusandhan National Research Foundation (ANRF), Government of India, through the ANRF Project CRG/2023/004921. 

\section*{Conflict of Interest}
The authors declare no conflicts of interest.

\noindent Komal Negi \\
Department of Mathematics, Indian Institute of Technology Ropar, Punjab, India.
\begin{verbatim} komal.20maz0004@iitrpr.ac.in \end{verbatim} 

\noindent Ayaka Shimizu \\ 
Osaka Central Advanced Mathematical Institute, Osaka Metropolitan University, Sugimoto, Osaka, 558-8585, Japan. 
\begin{verbatim} shimizu1984@gmail.com \end{verbatim}

\noindent Yoshiro Yaguchi \\ 
Maebashi Institute of Technology, Maebashi, Gunma, 371-0816, Japan. 
\begin{verbatim} y.yaguchi@maebashi-it.ac.jp \end{verbatim}

\noindent Madeti Prabhakar \\
Department of Mathematics, Indian Institute of Technology Ropar, Punjab, India.
\begin{verbatim} prabhakar@iitrpr.ac.in \end{verbatim}


\begin{thebibliography}{99}

\bibitem{BGRW} B. Baumeister, T. Gobet, K. Roberts, P. Wegener, On the Hurwitz action in finite Coxeter groups, \textit{J. Group Theory} {\bf 20} (2017), 103-131.

\bibitem{D} P. Dehornoy, Diagram colourings, \textit{Proceedings of the East Asian School of Knots, Links, and Related Topics February 16–20}, 2004.

\bibitem{LM} E. Liberman, M. Teicher, The Hurwitz Equivalence Problem is Undecidable, \textit{arXiv:math/0511153}, 2005.

\bibitem{KPS} K. Negi, M. Prabhakar, and S. Kamada, Twisted virtual braids and twisted links, \textit{Osaka J. Math.},  \textbf{61}(4), 2023.

\bibitem{KAM} K. Negi, A. Shimizu, M. Prabhakar, Warping labeling for twisted knots and twisted virtual braids, \textit{arXiv:2406.08505}, 2024.

\bibitem{KAY} K. Oshiro, A. Shimizu, Y. Yaguchi, Up-down colorings of virtual-link diagrams and the necessity of Reidemeister moves of type II, \textit{J. Knot Theory Ramifications}  {\bf 26} (2017), 1750073, 17 pp.

\bibitem{AY} A. Shimizu and Y. Yaguchi, Determinant of the OU matrix of a braid diagram, \textit{arXiv:2410.17778}, 2024(to appear in {\it J. Knot Theory Ramifications}).

\bibitem{Y} Y. Yaguchi, Determining the Hurwitz orbit of the standard generators of a braid group, \textit{Osaka J. Math.} {\bf 52} (2015), 59–70.

\end{thebibliography}
\end{document}